\documentclass[11pt]{article}
\usepackage{amsmath,amsfonts}
\newcommand{\VEC}{\vec}
\newcommand{\e}{{\VEC e}}

\newcommand{\ei}{\e_i}
\newcommand{\z}{\Vec z}
\newcommand{\lai}{\lambda_i}
\newcommand{\x}{\Vec x}
\newcommand{\y}{\Vec y}
\newcommand{\vu}{\Vec u}
\newcommand{\vr}{{\Vec r \, {}}}

\newcommand{\vv}{\Vec v}
\newcommand{\w}{\Vec w}

\newcommand{\field}[1]{\mathbb{#1}}
\def\bydef{\buildrel \triangle \over =}       

\def\epsilon{\varepsilon}
\def\R{\field{R}}
\def\P{\field{P}}

\newcommand{\BB}{{\cal B}}                    
\def\iot{I_{[0,T]}}

\def\iotr{\iot (\vr )}

\newtheorem{theorem}{Theorem}
\newtheorem{lemma}[theorem]{Lemma}
\newtheorem{cor}[theorem]{Corollary}

\newtheorem{assumption}{Assumption}

\numberwithin{equation}{section}

\newcommand{\blot}{\hfill{\vrule height .9ex width .8ex depth -.1ex }}
\newcommand{\EndPf}{\hfill $\blot$ \medskip}

\newcommand{\cali}{{\cal{I}}}

\newcommand{\va}{\Vec a}
\newcommand{\xstar}{\x^*}
\newcommand{\vq}{\Vec q}
\newcommand{\vrb}{\vr_B}
\newcommand{\vs}{\Vec s}
\newcommand{\ssb}{\vs_B}
\newcommand{\zn}{\z_n}
\newcommand{\vtheta}{\Vec \theta}
\begin{document}
\title{Uniqueness of a constrained variational problem and large
deviations of buffer size}
\author{Adam Shwartz, 
and
Alan Weiss 
\thanks{Adam Shwartz is the Julius M.\ and Bernice Naiman Chair in
Engineering, at the
Faculty of Electrical Engineering,
Technion---Israel Institute of Technology,
Haifa 32000, Israel. email adam@ee.technion.ac.il.
%
His work was supported in part by the fund for promotion of
research at the Technion and by the Dent charitable trust-non military research
fund.}
\thanks{Alan Weiss is with Mathworks, Inc., Natick, MA 01760-2098. email
aweiss@mathworks.com}}
%
%
%
\maketitle
\begin{abstract}
We show global uniqueness of the solution to a
class of constrained variational problems, using scaling
properties.
This is used to establish the essential uniqueness of solutions of a
large deviations problem in multiple dimensions.
The result is motivated by models of buffers,
and in particular the probability of, and typical path to overflow
in the limit of small buffers, which we analize.
\end{abstract}

\noindent
{\bf keywords}
Uniqueness of Variational problems, Large Deviations, Buffer overflow, AMS
model.
%
\section{Introduction}\label{s:intro}
We investigate
uniqueness of solutions to variational problems that arise in
sample-path large deviations. Our motivation comes from
models of buffers in telecommunication systems. The original
model was developed by Anick, Mitra, and Sondhi \cite{ams}.
Weiss \cite{We} cast this model in the
framework of sample-path large deviations
and showed for the probability of buffer overflow 
\begin{equation}\label{e:avarprob}
\lim_{n\to\infty} \frac 1 n \log \P \left (
  b_n (t) \ge B \right )
   =\inf_{(\vr, T)\in G(B)} \int_0^T
     \ell(\vr(t), \vr^\prime (t))\, dt ,
\end{equation}
where $\ell(\vr(t), \vr^\prime (t))$ is a positive ``cost" function, and
$G(B) $ is the set of paths $\vr $ ($d$-dimensional functions of $t$)
and terminal times $T$ satisfying a buffer overflow
property (see \S~\ref{s:buffer}).
Botvich and Duffield \cite{bd} (and independently Courcoubetis and Weber
\cite{cw} and Simonian and Guibert \cite{GS}) analyze
the probability of buffer overflow for more general
traffic models. For a superposition of $n$ traffic sources,
each bringing independent traffic $A(t)$ to the buffer
in time $(0, t)$, and with the buffer draining at rate $nc$ for some constant
$c$, they show that
the steady-state probability that the scaled buffer content
$b_n (t) = b(t) / n$ exceeds a level $B$ is approximately
\begin{equation}\label{e:bdresult}
\lim_{n\to\infty} \frac 1 n \log \P \left (
  b_n (t) \ge B \right )
     = - \inf_{t > 0} \sup_{\theta} \Bigl \{
     \theta (b + ct) - \log E\left (
       \exp(\theta A(t)) \right ) \Bigr \} .
\end{equation}
For more recent results see~e.g.~Mandjes and Mannersalo~\cite{MM} and
references therein.
Mandjes and Ridder~\cite{mrqs} showed how
to use the solution of \eqref{e:bdresult} to
find a minimal sample path $\vr(t)$ in \eqref{e:avarprob};
that is, a solution to the
variational problem. If this solution is unique, then this result
is quite useful, as shown by Freidlin and Wentzell \cite{fw}: if
$\z_n (t)$ represents the scaled state of the system at time $t$,
$b_n (t)$ the scaled buffer size, and
$\vr (t)$ is the {\em unique} solution to the variational problem
\eqref{e:avarprob}, then for each fixed $s \le T$,
\begin{equation}\label{e:spppty}
\lim_{n\to\infty}\P_{ss}\left ( \sup_{s\le t\le T} \left |
   \z_n (t) - \vr(t) \right | < \epsilon ~ \bigg | ~
    b_n (T) \ge B \right ) = 1.
\end{equation}
Here $\P_{ss}$ is steady-state probability
and $\P(A | B)$ is the probability of $A$ conditioned on $B$.
Thus, 
conditioned
on the buffer overflowing, we can determine beforehand
most likely how it did so.
So the sample-path approach yields interesting information when 
the minimizing $\vr(t)$ is unique.

The general problem of uniqueness for variational
problems is difficult. In general, it is
difficult to prove local uniqueness for variational problems,
and global uniqueness is even more difficult. In addition, our
problem is an optimization problem under a constraint (buffer exceeds $B$).
Constraints can often be handled using Lagrange multipliers, so that standard
uniqueness methods can be applied. However, this requires a proof
that the multiplier principle holds in this case. Note that due to
our application, the minimum must be shown to be unique in the
class of absolutely continuous paths, not simply $C^1$
or piecewise $C^1$ paths.
Finally, in terms of the application, there seem to be no direct
probabilistic methods for establishing uniqueness.

In this note we consider a constrained variational problem where
the constraint is given by a single functional.
We assume a certain scaling property of the
constraining functional together with homogeneity, convexity and
superlinear growth of the Lagrangian $\ell$ to show global
uniqueness of the minimizing solution. For our
motivating model this means the following. In previous studies of this model
(\cite{SWbook}, Chapter 13) we showed, through laborious
calculations, the applicability of the multiplier principle as
well as uniqueness, for a simple one-dimensional case. We show
here that the assumptions of our new uniqueness result hold for a
constant coefficient process. Then we approximate our model with a
constant coefficient one and show that, in the limit as the buffer size
$B\to 0$,
the variational problem has an {\em essentially\/} unique solution, in the
following sense. There is a unique path $\vr^* (t)$, $0\le t\le
T$, for the constant coefficient process,
such that any minimizing path $\vrb (t)$ of the real process
is close to a scaled version of $\vr^* (t)$:
\begin{equation}\label{e:firstresult}
\lim_{B\to 0} \frac 1{\sqrt{B}} \left ( \vrb (t\sqrt{B}) - \vrb(0) \right )
  + \vr^*(0) = \vr^*(t),\ 0\le t < T,
\end{equation}
uniformly in $t$, and $\vrb(0)\to\vr^* (0)$.
In other words, for each $\epsilon > 0$
there is a $\delta > 0$ such that for all $B < \delta$,
\begin{equation}\label{e:unifstate}
\sup_{0\le t\le T-\delta}\left | \frac 1{\sqrt{B}}
\left ( \vrb \left (t\sqrt{B} \right ) - \vrb(0) \right ) -
\left ( \vr^*(t) - \vr^*(0) \right ) \right | < \epsilon.
\end{equation}
The conditions under which this result holds are mild.

Note that this result together with a local uniqueness result for
the original model---for small enough $B$--- imply that there is a
unique solution for all $B$ small enough, for each given starting
point.
We do not give general conditions for the uniqueness of the
full buffer model, but some discussions appear in \S~\ref{s:weak}.
Our result does not imply full sample-path uniqueness
for the stochastic model with non-constant rates, 
although we are certain that uniqueness does hold. 
In paticular,
we have no effective bounds on the error (the relationship between
$\epsilon$ and $\delta$ in \eqref{e:unifstate}).

In \S~\ref{s:cc} we introduce the general variational problem, and
in~\S~\ref{s:assns} we state the abstract assumptions under which uniqueness
holds. Sections~\ref{s:yesitexists} and~\ref{s:onlyone} prove existence and
uniqueness resp. In \S~\ref{s:constraint} we derive a property of the shape of
optimal paths for buffer problems. In \S~\ref{s:buffer} we formulate the
stochastic
buffer problem and obtain some properties, in~\S~\ref{s:weak} we prove weak
uniqueness and in~\S~\ref{s:strong} we show a strong uniqueness result and
discuss the desired full uniqueness result.

\section{The variational problem}\label{s:cc}
We assume throughout that the functional
$\iotr$
on absolutely continuous paths $\vr(t)$
($d$-dimensional functions
of a real variable)
has a representation
\begin{equation}\label{e:iotrform}
\iotr \bydef \int_0^T \ell\left ( \vr(t), \frac d{dt}\vr(t)\right )\, dt,
\end{equation}
where $\ell(\x, \y)$ is a positive function, convex in $\y$.
The {\em constant coefficient\/} problem has
$\ell(\x, \y) = \ell(\y)$---independent of $\x$.
We prove our uniqueness result for a general class of
constraints. So, let $ \BB $ be a given positive functional
on absolutely continuous paths.
Our optimization problem is 
\vspace{3 mm}

\noindent{\bf Problem A.} Given $\x$ and $B > 0$, find a time $T$ and
an absolutely continuous path $\vr(t)$ that minimize
$\iotr$ subject to $\BB(\vr, T) \ge B$.
\vspace{3 mm}

\noindent
Denote this minimal cost by ${\cal I\/}$.
This is a ``free time" problem; hence
straightofrward convexity arguements cannot be used to infer uniqueness.

\section{Assumptions}\label{s:assns}
Throughout, we use the $ \sup $ norm on the space of (measurable) functions.
The following assumptions apply to the constant coefficient problem.
We shall use all 
save the last assumption in
order to prove existence of solutions to Problem A,
and a different subset to prove uniqueness of solutions. 
\begin{assumption}\label{a:minasone}
\begin{itemize}
\item[]
\item [{[A]}] $\ell(\y)$ is positive and strictly convex in $\y$.
\item [{[B]}] ${\displaystyle \lim_{|\y|\to\infty}
     {\ell(\y)}\slash{|\y|} = \infty}$.
\item [{[C]}] $\BB(\vr, t)$ is positive, continuous in
$ (\vr , t )$, and $\BB(\vr, 0) = 0$.
\item [{[D]}]
For some $\vs$,
$\ell (\vs) < \infty$
and
$\lim_{t\to\infty}\BB(\vr, t) = \infty$
          for $\vr(t) = \x + t\vs$.
\item [{[E]}] 
   Given $\x$ and $B$ there is a 
   $T_1 < \infty$
   such that any path $\vr(t)$ 
   on $0\le t\le T$ with $T > T_1$, 
   $\vr(0) = \x$ and $\BB(\vr, T) \ge B$
   has $\iotr > {\cal I} + 1$.
\item [{[F]}] For all $\vr(t)$ and $\z(t)$ on $0\le t\le T$
with $\vr(0) = \z(0)= \x$, 
$ \BB $ satisfies
\begin{enumerate}
\item For any $\gamma\in(0, 1)$ and $0\le t\le T$ we have
\begin{equation}\label{e:Bconcave}
\BB(\gamma\vr + (1-\gamma)\z, t) \ge \gamma\BB(\vr, t) + (1-\gamma)\BB(\z, t).
\end{equation}
\item There exists a positive increasing function $f(\alpha)$, defined
for $\alpha > 0$,
that satisfies $f(\alpha)f(1/\alpha) \equiv 1$, 
and is either the identity ($f(\alpha) = \alpha$), or is
strictly convex for $\alpha \ge 1$
such that for every $\vr(t)$ on $0\le t\le T$, the scaled path
$\y(t)\bydef\alpha\vr(t/\alpha)$
satisfies
\begin{equation}\label{e:scalinglaw}
\BB(\y, t) = f(\alpha) \BB(\vr, \alpha t) \quad
0\le t\le \alpha T
\end{equation}
\end{enumerate}
\end{itemize}
\end{assumption}
Note that Assumption \ref{a:minasone} [A] implies
for any $\gamma\in(0, 1)$ we have
\begin{equation}\label{e:convexi}
\iot (\gamma\vr + (1-\gamma)\z)\le\gamma\iotr + (1-\gamma)\iot(\z),
\end{equation}
where the inequality is strict unless $\vr(t) \equiv \z(t)$.
Furthermore, the representation of $I$ as an integral implies that,
for $\y$ defined in terms of $\vr$ as in Assumption \ref{a:minasone} [F] 2,
$ I_{[0, \alpha T]}(\y) = \alpha\iotr $.

The prototypical $f(\alpha)$ satisfying 
[F] 2 is $f(\alpha) = \alpha^\beta$ for $\beta > 1$.
Our motivating example has $\beta = 2$; see \eqref{e:genBBdef}.
The conditions 
do not imply that $f(\alpha)$ is convex for $\alpha < 1$.
They imply that $f(1) = 1$ and $f(\alpha)\to 0$ as $\alpha\to 0$.

\section{Existence}\label{s:yesitexists}
%
The real import of this paper
is uniqueness; for completeness, we prove that
for the constant coefficient problem, there exist
solutions to Problem A under our assumptions.

\begin{theorem}\label{bigdeal}
Under Assumption \ref{a:minasone} parts [A]--[E],
there exists a solution to Problem A.
\end{theorem}

{\bf Proof}:
By Assumption \ref{a:minasone} [D]
the linear function $\vr(t) = \x + t\vs $
has $\iotr < \infty$
and makes $\BB(\vr, T) \ge B$ for $T$ large enough.
Therefore we are not
minimizing over an empty set.
By Assumption \ref{a:minasone} [E]
we may restrict to minimizations over sets of bounded time.

Take a minimizing sequence $\vr_i$ on $[0, T_i]$
with $\BB(\vr_i, T_i) \ge B$ and $\vr_i (0) = \x$.
To prove existence of an optimal $\vr$ we need to
show that a converging subsequence of approximate minimizers has a
minimal limit. By \cite[Lemma 5.18]{SWbook}, under
Assumption \ref{a:minasone} [A] and [B], a set of paths $\vr_i (t)$
over a bounded interval $[0, T]$,
having $\vr_i (0)$ in a compact set, and with
bounded $I$-functions, is uniformly absolutely continuous, hence pre-compact.
Therefore, since any minimizing sequence has bounded $T$,
there exists a convergent subsequence, and the limiting path
$\vr(t)$ is absolutely continuous.
Again under Assumption \ref{a:minasone} [A] and [B],
\cite[Lemma 5.42]{SWbook} shows that the functional $\iot$
is lower semicontinuous in $\vr$; therefore, the limiting path satisfies
$\iotr \le \lim_i\iot(\vr_i)$.
The continuity of $\BB$ in $(\vr , T)$ shows that
$\BB(\vr, T) \ge B$.
\EndPf

\section{Uniqueness}\label{s:onlyone}
In this section we establish that the solution to
the constant coefficient problem 
is unique.
Note that this is only interesting because Problem A has the constraint
$\BB(\vr, T) \ge B$ and $T$ is not determined.
If we were simply trying to minimize $\iotr$ with fixed initial
and final points, it is trivial to show that the unique solution is
a straight line path (e.g., \cite[Lemma 5.15]{SWbook}).

\begin{theorem}\label{t:itsunique}
Fix $\x$ and $B > 0$.
Suppose Assumption \ref{a:minasone} [A], [B], [C], and [F] hold.
If $\vr(t)$ $0\le t\le T$ and $\y(t)$ $0\le t\le T_1$
solve Problem A,
then $T = T_1$ and $\vr(t) = \y(t)$, $0\le t\le T$.
\end{theorem}
Note: in the proof, we need to rescale time in part because the
optimal paths may live on different time intervals.

{\bf Proof}:
Without loss of generality assume $T_1 \geq T$.
Define $ \alpha = T_1 \slash T \geq 1$.
By assumption,
\begin{align}\label{e:vrvyequal}\notag
\iotr &= \cali & \BB(\vr, T) &\ge B \\
\notag
I_{[0, \alpha T]}(\y) &= \cali & \BB(\y, \alpha T) &\ge B .
\end{align}
Take
$
\vu (t) \bydef \alpha \vr(t/\alpha).
$
Then 
$
I_{[0, \alpha T]} (\vu) = \alpha \cali  
$ and $
\BB(\vu, \alpha T) \ge f(\alpha) B .
$
Now for any $0<\gamma<1$ define
$
\vv = \gamma \y(t) + (1-\gamma) \vu (t).
$
Then we have
\begin{equation}\label{e:vvumore}
I_{[0, \alpha T]} (\vv) \le (\gamma + (1-\gamma)\alpha) \cali , \quad
\BB(\vv, \alpha T) \ge (\gamma + (1-\gamma)f(\alpha)) B .
\end{equation}
Now scale $\vv$ to get 
$
\w(t) = \delta \vv(t/\delta).
$
Then we have
\begin{equation}\label{e:vvwscale}
I_{[0, \delta\alpha T]}(\w) \le \delta(\gamma + (1-\gamma)\alpha) \cali ,\quad
\BB(\w, \delta\alpha T)\ge f(\delta)(\gamma + (1-\gamma)f(\alpha)) B
\end{equation}
where the first inequality is strict unless $\vu = \y$.
Note that since $f(1) = 1$, so $\gamma + (1-\gamma)f(\alpha) =
\gamma f(1) + (1-\gamma)f(\alpha)$, we have by \eqref{e:vvwscale} that
\begin{equation}\label{e:newscale}
\BB(\w, \delta\alpha T) > f(\delta)f(\gamma + (1-\gamma)\alpha) B
\end{equation}
by the strict convexity of $f$, as long as $\alpha > 1$ and $f$
is not the identity.
If there is strict inequality in \eqref{e:newscale}, we may
find a $\delta < 1/(\gamma + (1-\gamma)\alpha)$ so that
$\BB(\w, \delta\alpha T) > B$.
But then by \eqref{e:vvwscale}, $I_{[0, \delta\alpha T]}(\w) < \cali$,
contradicting the assumed minimality of $\cali$.
If 
$f$ is the identity, so that inequality \eqref{e:newscale}
is not strict, then we choose $\delta = 1/(\gamma + (1-\gamma)\alpha)$,
and get a contradiction to the assumed minimality of $\cali$
from the strict inequality in \eqref{e:vvwscale}.
\EndPf

\section{Integral constraint: concave buffer}\label{s:constraint}
This section
illustrates our assumptions on $\BB$, and derives a result of independent
interest for the study of properties of buffer models.
This result is not needed for the main argument.
We consider
the constant coefficient problem,
and a functional $\BB$ that is a
the traditional buffer model,
and show that optimizing paths have to be concave in a preferred direction.
For a given vector $\va$ with positive components, we take
\begin{align}\label{e:genBBdef}
\frac{d}{dt} \BB(\vr, t) &= \begin{cases}
\langle\vr(t), \va\rangle &\text{if }\langle\vr(t), \va\rangle > 0
\text{ or }\BB(\vr, t) > 0\\
0 &\text{otherwise},
\end{cases}\\
\label{e:Bstartsempty}
\BB(\vr, 0) &= 0.
\end{align}
Clearly $\BB$
satisfies Assumption~\ref{a:minasone}
[C]--[D]\footnote{[D] holds
provided $\ell(\vs)<\infty$ for some $\vs$ with
$\langle\vs,\va\rangle > 0$.};
we discuss [E] later.
We now show [F].
\begin{lemma}\label{l:triviality}
For any integrable function $g(t)$, we have
\begin{equation}\label{e:gscale}
u_\alpha (t) \bydef \int_{\alpha s}^{\alpha t} \alpha g(w/\alpha)\, dw
= \alpha^2 u_1 (t).
\end{equation}
\end{lemma}
\begin{lemma}\label{l:bisanintegral}
$\BB(\vr, t)$ defined by \eqref{e:genBBdef}, \eqref{e:Bstartsempty} satisfies
\begin{equation}\label{e:alternateBB}
\BB(\vr, t) = \sup_{0\le s\le t}\int_s^t \langle\vr(w), \va\rangle\, dw.
\end{equation}
\end{lemma}
\begin{cor}\label{c:Bscales}
Given $\vr(t)$ set
$\y(t)\bydef\alpha\vr(t/\alpha)$. Then
$\BB(\y, T) = \alpha^2 \BB(\vr, \alpha T)$.
\end{cor}
These results are immediate;
for a proof of Lemma \ref{l:bisanintegral} see~\cite[\S
13.2]{SWbook}, \cite{SWtimeScale}.

We begin with an intuitive result: 
optimal paths 
never have $\langle\vr(t), \va\rangle < 0$.
If they did, 
then the buffer
could be increased and the cost $I$ decreased by eliminating the point.
The proof is exactly this argument, but we must account for the
shifting of time and position
(shift segments of paths to make them continuous)
that come from eliminating intervals of time.
\begin{lemma}\label{l:bufferpositive}
Any minimal $\vr(t)$, $0\le t\le T$ with
$\langle \vr(0), \va\rangle \ge 0$ satisfies
$\langle \vr(t), \va\rangle \ge 0$, $t\in[0, T]$.
\end{lemma}

{\bf Proof}:
Eliminate any times with $\langle \vr(t), \va\rangle < 0$, then
reconstitute $\vr$. 
Mathematically, let
\begin{equation}\label{e:changeotime}
s(t) = \int_0^t 1_{\langle \vr(u), \va\rangle \ge 0} \, du ;
\end{equation}
Let
$
\tau(q) = \inf\{t~:~s(t) = q\}
$
be the 
first time by which 
$\langle \vr(u), \va\rangle \ge 0$
for a total time of $q$. Then
$ \tau (s (p)) = s ( \tau (p) ) = p $, that is, $\tau $ and $s$ are inverses.
Then define
\begin{equation}\label{e:ronenew}
\vr_1 (t) = \vr(0) + \int_0^{\tau(t)}\frac{d}{du}\vr(s(u))\, du.
\end{equation}
This is well defined since $\vr $ is absolutely continuous.
The path $\vr_1 (t)$ has the same increments as the path $\vr(t)$
once the times $\langle \vr(t), \va\rangle < 0$ are eliminated,
and it is continuous.
Clearly $\iotr \geq I_{[0, s(T)]} (\vr_1 )$.
But if $\BB(\vr, T) = B> 0$ then $\BB(\vr_1, s(T)) \geq B$,
since the two functions are equal at corresponding (shifted) times
whenever $\langle\vr(t),\va\rangle\ge 0$;
see~\eqref{e:genBBdef} or~\eqref{e:alternateBB}.
In fact, if $\langle\vr(t),\va\rangle < 0$ anywhere, then
$ \BB(\vr_1, s(T)) > B$. But in this case, by the continuity
of $ \BB $ (assumption~\ref{a:minasone} [C]) and its definition, there is some
$T_1 < s(T) $ so that $ \BB(\vr_1, T_1 ) = B$ and, by
definition~\eqref{e:iotrform}, since $\ell $ is positive,
$ \iotr > I_{[0, T_1 ]} (\vr_1 )$, so that $\vr $ is not optimal
.
\EndPf

\begin{lemma}\label{l:concavity}
Given $B$, $T$ and $\x$ with $\langle\x,\va\rangle\ge 0$,
any $\vr (t)$ with $\vr(0) = \x$ that has $\BB(\vr, T) = B$ and has
minimal cost $\iotr$, has the property that the function $\langle
\vr (t), \va\rangle$ is concave, $0\le t\le T$.
\end{lemma}

{\bf Proof}: Suppose that the function $\langle \vr (t), \va\rangle$ lies
strictly below its concave envelope for some $L < t < M$.
Let $v(t)$ denote the concave envelope of $\langle \vr (t),
\va\rangle$; hence $v(t)$ is linear on $[L, M]$.
By~\cite[Lemma 5.15]{SWbook}, replacing $ \vr $ with a linear
function over the interval $ [L,M ] $ strictly decreases $ I $, and
since $v(t) > \langle\vr(t),\va\rangle$ for $L < t < M$,
by \eqref{e:genBBdef}, $\BB(\vr, t)$ increases for $t>L$ .
\EndPf

\section{Stochastic buffer model}\label{s:buffer}
This section gives an overview of the connection between
a class of stochastic traffic and buffer models
and the variational Problem A. Precise definitions and
assumptions will be given later.
There are two components in the models:
a Markov chain $\x(t)$ (called the state of the system)
and the buffer itself.
The Markov chain operates in continuous time
and its state lies in the $K$-dimensional lattice
with positive integer components.
There are $J$ possible transitions, with
transition rates $\lambda_i (\x)$
in direction $\ei$, $1\le i\le J$.
This means that the state jumps from $\x$ to $\x + \ei$ with
Poisson rate $\lai (\x)$.
The interpretation of the state
is that there are $x_i$ sources of traffic in state $i$.
There is a single buffer with continuous contents $B(t)$.
The buffer contents are filled by each source in state $i$
at rate $a_i >0$, and the total buffer drain rate is $C$.
Mathematically, the buffer satisfies the following equation
(cf.~\eqref{e:genBBdef}--\eqref{e:alternateBB}):
\begin{equation}\label{e:Bequation}
\frac d{dt} B(t) = \begin{cases}
\langle \x (t), \va\rangle - C &
         \text{if }\langle \x (t), \va\rangle > C \text{ or }B(t) > 0\\
0 & \text{otherwise.}
\end{cases}
\end{equation}
If $B(0) = 0$ then it is easy to see that
\begin{equation}\label{e:BBaltdef}
B(t) = \sup_{0\le s\le t} \int_{u=s}^t \langle a, \x (u)\rangle \, du .
\end{equation}
We define $\BB(\x, t)$ as this $B(t)$ for a given (state) sample path $\x(t)$.

As an example, suppose that there are $K=2$ types of traffic sources,
phone and data.
Type $1$ sources, phone, generate traffic at rate $1$, arrive at
Poisson rate $\lambda$, and depart at rate $\mu$ each.
Type $2$ sources, data, generate traffic with rate $5$,
arrive at rate $\theta$, and depart at rate $\psi$ each.
Furthermore, a phone source can turn into a data source
(a person talks for a while, then starts sending a fax);
this occurs at rate $\gamma$ for each source in state $1$,
and, of course, a data source can become a phone source
with rate $\delta$.
We suppose the buffer drain rate is 100.
This leads to a model with the following jump rates and transitions:
\begin{align*}
\e_1 &= (1, 0) & \lambda_1 &= \lambda &
\e_2 &= (-1, 0) & \lambda_2 &= \mu x_1\\
\e_3 &= (0, 1) & \lambda_3  &= \theta  &
\e_4 &= (0, -1) & \lambda_4  &= \psi x_2  \\
\e_5 &= (-1, 1) & \lambda_5  &= \gamma x_1  &
\e_6 &= (1, -1) & \lambda_6  &= \delta x_2
\end{align*}
\begin{equation*}
\frac{d}{dt} B(t) = \begin{cases}
x_1 (t) + 5 x_2 (t) - 100 & \text{if }x_1 (t) + 5 x_2 (t) > 100
   \text{ or }B(t) > 0\\
0 & \text{otherwise.}
\end{cases}
\end{equation*}
This is an open model, since sources arrive and depart.
For $\lambda_i = 0$, $1\le i\le 4$, this is a closed model:
transitions
are between types of sources, and the total number of sources
is constant.
We could also allow for correlated arrivals and
departures. For example, we could have a person arrive and
start to talk on the phone and use his email at the same time;
this would be an arrival with $e_7 = (1, 1)$.
Similarly, this type of person could depart $(e_8 = (-1, -1))$,
stop one of his activities ($e_2$ or $e_4$),
or a couple could start their day together with $e_9 = (2, 2)$ etc.

We suppose that the process $\x(t)$ may be scaled to $\zn(t)$
with jumps of $\zn$ having size $\ei / n$ with rate
$n\lai (\zn (t))$.
The scaled buffer $b_n (t)$ is defined similarly to \eqref{e:Bequation}:
\begin{equation}\label{e:scaledBeq}
\frac d{dt} b_n(t) = \begin{cases}
\langle z_n (t), \va\rangle - C &
         \text{if }\langle z_n (t), \va\rangle > C \text{ or }b_n(t) > 0\\
0 & \text{otherwise.}
\end{cases}
\end{equation}
In some models $C$ is also scaled by $n$
so as to have a finite limit as $n\to\infty$.
In general we suppose that the jump rates $\lai(\x)$ keep each component
of $\zn(t)$ positive; that is, we do not consider ``sources" that
might drain the buffer.
Clearly $b_n (t) = \BB(\zn, t)$ under appropriate scaling.

We showed in \cite{diminish} that, under a wide variety
of assumptions on the $\lambda_i(\x)$, the process
$\zn (t)$ satisfies a large deviations principle.
Furthermore, we showed in \cite[Chapter 13]{SWtimeScale,SWbook}, that
in a buffer model very similar to the one just described,
statistical properties of the buffer size process may also be approximated by
solutions to variational problems.
(There have been many other analyses of these problems;
typically, they consider steady state statistics. See, e.g.\
Mandjes and Ridder \cite{mrqs} or, for a non-sample-path approach
Botvich and Duffield \cite{bd}.)
In particular, it was shown that
the steady-state probability that the buffer
exceeds a level $B$ is given by $\exp(-nI + o(n))$, where
$I$ is the solution to a variational problem we describe below.

\subsection{The variational problem for Buffer overflow}
Given the rates and jump directions,
define the local cost function by
\begin{equation}\label{e:elldef}
\ell(\x, \y) \bydef \sup_{\vtheta\in\R^K}\left ( \langle\vtheta, \y\rangle
   - \sum_i \lambda_i (\x) \left ( \exp\langle \vtheta, \ei\rangle -1 \right ) \right ).
\end{equation}
Then define the action functional, or cost 
of an absolutely continuous function $\vr (t)\in\R^K$ 
\begin{equation}\label{e:Idef}
\iotr \bydef \int_0^T \ell\left ( \vr(t), \frac d{dt} \vr(t) \right )\, dt.
\end{equation}

The variational problem associated with making a buffer exceed level $B$
is given as follows. We assume that the process is typically near its
steady state value, denoted by $ \vq $, and that at that point the capacity
$C$ suffices to serve the traffic, namely $ \langle \vq, \va\rangle < C$.
Let ${\cal I}$ be defined by
\begin{equation}\label{e:calIdef}
{\cal I} = \inf\left\{\iotr~:~\vr(0) = \vq,\ T < \infty,\ \BB(\vr, T) \ge B\right\}.
\end{equation}
Usually the optimizing $T$ is infinite, since the time
to leave $\vq$ is infinite, so we
perform a standard shifting of time so that the upcrossing
of the path $\vr$ to the hyperplane $\langle\vr, \va\rangle = C$ takes
place at time $0$. The time $T$ after time $0$ required to fill the buffer
is then typically finite.
If there is a unique solution $(T, \vr)$ with $I_{[-\infty, T]} (\vr)= {\cal I}$
and $\BB(\vr, T)\ge B$ then equation \eqref{e:spppty} holds; in other words,
given a buffer overflow, we know how it was almost certain to have happened.
For a proof of this assertion see \cite[Chapter 13]{SWbook}.

The problem of finding the minimal cost path $\vr$ starting from $\vq$ and making
$\BB(\vr, T) \ge B$ naturally splits into finding the cheapest path
from $\vq$ to sone $\x$ with $\langle \x, \va\rangle = C$, followed by the
cheapest path starting from $\x$ making $\BB(\vr, T) \ge B$, and then minimizing
the total cost $I$ over all such $\x$.
This is clearly valid because the buffer cannot begin to increase until
$\langle\vr,\va\rangle > C$.
It is also clear that the cost $\iotr$ of making a buffer $B$
tends to zero as $B\to 0$ (since we can take a small period of time
to make the small buffer, and $\ell$ is finite for some directions).
Therefore, if there is a unique upcrossing point $\xstar$, then for small
$B$, the optimal point $x$ (in the full problem of making a small buffer)
will be close to $\xstar$, the optimal upcrossing point.

\subsection{Typical behavior}
To describe the typical behavior of the state 
and buffer size processes,
define 
the process $\z_\infty (t)$
\begin{equation}\label{e:vvxdef}
\vv(\x) \bydef \sum_i \lambda_i(\x) \ei ,
\quad
\frac{d}{dt}\z_\infty (t) = \vv(\z_\infty (t) ).
\end{equation}
Kurtz's theorem
\cite[Chapter 5.1]{SWbook}
states that, over finite time intervals,
the process $\zn(t)$ is extremely likely to remain close to
the path $\z_\infty (t)$ that starts at the same point $\z_\infty (0) = \zn(0)$.
The intuition is that $\vv $
is the average drift of the process $\zn (t)$,
in the sense that
\begin{equation}\label{e:zndrifter}
\lim_{\delta\downarrow 0}\lim_{n\to\infty}
\frac{\zn(t + \delta) - \zn(t)}{\delta} - \vv(\zn (t)) = 0
\text{ with probability one.}
\end{equation}
%
From this one can infer~\cite[Cor.\ 8]{SWtimeScale} that the buffer size
process $ b_n $ is also extremely likely to remain close to a limiting
path, which can be calculated.

We will show next that given an initial point $\x$ with
$\langle \x, \va\rangle = C$, and given a small $B$,
there is an {\em essentially unique}
path $\vr$ with $\vr(0) = \x$ and $\BB(\vr, T)\ge B)$.
This clearly does not solve the full uniqueness problem, because
although we know that $\x\approx\xstar$, for each $B$
we have no idea whether there is a unique associated point $\x$.

\section{Weak uniqueness}\label{s:weak}
We now state the assumptions under which the probabilistic problem
of the buffer being large is described by the deterministic
variational Problem A.
\begin{assumption}\label{a:probassns}
\begin{enumerate} \item[]
\item[A.] The processes $\zn(t)$ and $ b_n (t) $ satisfy a sample-path
large deviations principle;
the rate function for $ \zn $ is $\iot$.
\item[B.] There is a unique globally attracting point $\vq$ of $\z_\infty (t)$.
\item[C.] $\langle\vq, \va\rangle < C$.
\item[D.] The Freidlin-Wentzell theory~\cite[Chapter 6]{SWbook} applies to
both $\zn (t)$ and 
$b_n (t)$.
\item[E.] There is a unique 
endpoint
$\xstar$ 
of the upcrossing problem from $\vq$ to 
$\{ \x : \langle \x, \va\rangle = C \}$.
\item[F.] $\langle \xstar, \vv(\xstar)\rangle < 0$.
\end{enumerate}
\end{assumption}
Since $ \BB $ is a continuous functional, the assumption that the buffer
process satisfies a large deviations principle is a consequence of
the fact that $ \z_n $ satisfies the principle, and the rate function for
$ b_n $ can be derived from $\iot $~\cite[\S 2.3]{SWbook}:
see e.g.~\cite[Cor.\ 8]{SWtimeScale}.

\begin{theorem}\label{t:FWapplies}
Let $ \P_{ss} $ denote steady state probability.
Under Assumption \ref{a:probassns},
\begin{equation}\notag 
\lim_{n\to\infty}\frac 1 n \log \P_{ss}(\BB(\zn, t)\ge B)
  = -\inf_{\vr , T}
      \left \{ \iotr : \vr(0) = \vq,\ \BB(\vr, T)\ge B \right \} .
\end{equation}
Furthermore, if $\vr$ is unique to within a time shift, then
(shifting time so that the buffer overflow takes place at time $t = 0$)
the estimate~\eqref{e:spppty} holds.
That is, $\zn (t) \approx \vr(t)$
over each fixed interval of time before the time of overflow,
assuming that overflow occurs at time $T$.
\end{theorem}

The proof of Theorem \ref{t:FWapplies} is standard from the
Freidlin-Wentzell theory, and will not be given here;
see \cite{fw} or \cite[Chapter 6, 13]{SWbook}.
However, we must explain when Assumption \ref{a:probassns} can hold.
Part A was proved in \cite{diminish}, under the following conditions.
The process $\zn (t)$ is assumed to live on the positive quadrant $G$
(in fact, the domain $G$ in \cite{diminish} is quite
general). 
The jump rates $\lai(\x)$ were assumed Lipschitz continuous,
strictly positive in the interior of the domain $G$,
and $\lai(\x)\to 0$ as $x\to\y\in\partial G$ for every $i$
such that a jump in direction $\ei$ from $\y$ would take
$\zn$ out of $G$.
The associated $\ei$ were assumed to positively span all of $\R^K$
(that is, for any $\x\in\R^K$ there are $\beta_i \ge 0$
such that $\x = \sum_i \ei\beta_i$).
There were some further technical conditions that are quite minor
and will not be described here.
Part B of Assumption \ref{a:probassns} is simply that
the process $\zn (t)$ tends to go to a small neighborhood of $\vq$.
Part C means that the buffer $B(t)$ tends to decrease when the process is
in its most likely state.
Part D implies that steady-state quantities can be calculated from
the transient sample-path large deviations
principle~\cite[Chapters 6,13]{SWbook}.
Part F means once the process $\zn(t)$ reaches the most likely
point where the buffer can begin to fill, the process tends to
go in a direction that has the buffer decrease.
This means that, most likely, even when the buffer begins to fill,
it will almost immediately begin to decrease.

The 
assumptions can be largely eliminated if we change the question we ask.
Suppose we want to know how the buffer fills, starting from
some point $\x$ with $\langle\x,\va\rangle = C$, and how likely
it is to fill to a level $B$.
\begin{theorem}\label{t:transt}
Under Assumption \ref{a:probassns} A,
given $\x$ with $\langle\x,\va\rangle = C$,
\begin{equation}\notag 
\lim_{n\to\infty}\frac 1 n \log \P(\BB(\zn, T)\ge B)
  = -\inf\left \{ \iotr~:~\vr(0) = \x,\ T,\ \BB(\vr, T)\ge B\right \}.
\end{equation}
Furthermore, if $\vr$ (and $T < \infty$) are unique
then for every $\epsilon > 0$ 
\begin{equation}\label{e:uniquebufferpath}
\lim_{n\to\infty}\P\left (\sup_{0\le t\le T} \left | \zn(t)
- \vr(t) \right | < \epsilon ~\bigg | ~\BB(\zn, T)\ge B\right ) = 1.
\end{equation}
\end{theorem}
The proof is standard: see~\cite{SWtimeScale}.\\
It is because of Theorem \ref{t:FWapplies} and~\eqref{e:uniquebufferpath} of
Theorem~\ref{t:transt} that we are motivated to
study uniqueness for solutions to Problem A.
This concludes our discussion of the connection between the probability
and variational problems.
We return to the study of the variational problem.

Suppose that there is a unique upcrossing point $\xstar$, the lowest cost
point from $\vq$ to the hyperplane $\langle\x,\va\rangle = C$.
Consider the local cost function $\ell(\xstar, \y)$, which is defined
for any $(\x, \y)$ in \eqref{e:elldef}.
When viewed as a constant coefficient cost $\ell(\y)$,
by~\cite[Chapter 5.2]{SWbook} $\ell(\xstar, \y)$ 
satisfies Assumption~\ref{a:minasone} [A], [B] and [E] and, for nontrivial
Markov models, also [D].
When combined with the functional $\BB$ given in \eqref{e:Bequation},
it is easy to check that all of Assumption \ref{a:minasone} [A]--[F] hold
(see \cite[Chapter 5.2]{SWbook} for properties of $\ell(\y)$,
Chapter 13 for properties of $\BB$).
Therefore, the solution to Problem A for constant coefficient cost function,
with paths starting at $\xstar$
making buffer at least $B$, exists and is unique.
For the case $B = 1$ we call the resulting solution $\vr(t)$, and the time $T$.
Actually, the path may start at any $\x_1$ with $\langle\x, \va\rangle = C$,
as long as we take the cost function given by $\ell(\xstar, \y)$, not $\ell(\x_1, \y)$,
and then the solution will be $\vr(t) + \x_1 - \xstar$, since $ \iot $ is
invariant under shifts, while $ \BB $ is invariant under shifts along
the hyperplane $\{ \x : \langle \x , \va \rangle = 0 \} $.
Furthermore, our scaling property gives the optimal path for any $B\neq 1$
as a scaled version of $\vr$, namely
$\alpha\vr(t/\alpha) + \xstar - \alpha\xstar$, for $0\le t\le T\alpha$,
where $\alpha = \sqrt{B}$.

We now switch back to cost function $\ell(\x, \y)$ depending on $\x$.
Many properties of this function are detailed in \cite[Chapter 5]{SWbook}.
If the positive cone spanned by the $\ei$ is the entire space $\R^K$,
and if all $\lai (\x) > 0$, then $0\le\ell(\x, \y) < \infty$, and $\ell(\x, \y)$
is convex in $\y$.
\cite[Thm.\ 5.35]{SWbook} shows that for any $\epsilon > 0$, $c > 0$,
there exists a $\delta > 0$ such that if $J_{[0, T]} (\vr)$ is defined to be the same
as $\iotr$ except using the local cost $\ell(\xstar, \frac d{dt} \vr(t) )$,
and if $J_{[0, T]} (\vr) \le c$ and $|\lai (\vr(t)) - \lai(\xstar)| < \delta$ then
$|\iotr - J_{[0, T]} (\vr)| < \epsilon$. In other words, the functional $\iot$
is uniformly continuous in sets of bounded cost.
Also, \cite[Prop.\ 5.46]{SWbook} shows that for any $c, T < \infty$, the set of
$\vr (t)$ with $\vr (0)$ in a compact set, and with $\iotr \le c$, is compact.
Compactness is in the sup norm topology for the paths in $[0, T]\to\R^K$.

The basis of our argument is that when $B$ is small,
paths $\vrb(t)$, $0\le t\le T_B$,
can stay very near $\vrb(0)$ and have $T_B$
small, and still make $\BB(\vrb, T_B) = B$.
So the local cost function $\ell(\vrb(t), \frac{d}{dt}\vrb(t))$ is
very nearly equal to $\ell(\vrb(0), \frac{d}{dt}\vrb(t))$,
a constant coefficient cost, for which we know uniqueness applies.
This motivates the following cost function, which is based on
considering a small neighborhood of $\vrb(0)$.
Given $B$ and $\vrb(0)$ with $\langle\vrb(0), \va\rangle = C$, define
\begin{equation}\label{e:lbdef}
\ell_B (\x, \y) \bydef \ell\left (\vrb(0) + \sqrt{B}(\x - \vrb(0)), \y\right ).
\end{equation}
For any such path $\vrb (t)$ 
that has $\BB (\vr(T_B)) = B$ we define a scaled path
\begin{equation}\label{e:ssbdeff}
\ssb (t) \bydef \vrb (0) + \frac{1}{\sqrt{B}} \left (\vrb(t\sqrt{B}) - \vrb (0) \right),
\ 0\le t\le T_B / \sqrt{B} .
\end{equation}
A simple calculation shows that
\begin{equation}\label{e:ssbbdd}
\BB \left(\ssb, T_B / \sqrt{B}\right ) = 1.
\end{equation}
Now we find
\begin{multline}\label{e:ellBcalc}
\int_0^{T_B / \sqrt{B}} \ell_B \left (\ssb(t), \frac{d}{dt} \ssb (t)\right )\, dt \\
\begin{split}
&= \int_0^{T_B / \sqrt{B}} \ell_B \left (\vrb (0) + \frac{1}{\sqrt{B}}
     (\vrb(t\sqrt{B}) - \vrb (0) ), \frac d{dt} \vrb (t\sqrt{B})\right )\, dt\\
&= \int_0^{T_B / \sqrt{B}} \ell\left (\vrb (0) +  (\vrb(t\sqrt{B}) - \vrb (0) ),
     \frac d{dt} \vrb (t\sqrt{B})\right )\, dt\\
&= \frac{1}{\sqrt{B}} \int_0^{T_B} \ell\left (\vrb(u), \frac d{du} \vrb (u)\right )\, du.
\end{split}
\end{multline}

It is easy to show that, under Assumption \ref{a:probassns} E,
as $B\to 0$, minimal solutions to
Problem A starting from $\vrb(0)$ with $\langle\vrb(0), \va\rangle = C$
must have $T_B/\sqrt{B}$ bounded, since there are paths $\vr (t)$
that have $T/\sqrt{B}$ bounded, and have cost of order $\sqrt{B}$, and
any path that stays above $\langle\vr(t), \va\rangle = C$ for time $T$
has cost that is bounded below by a function linear in $T$.
Therefore, if we consider solutions to Problem A with buffer size $B$ small,
then the scaled problem with $\ell_B$ cost function
has solution $\ssb$ in bounded time $T_B/\sqrt{B}$.

\begin{theorem}\label{t:samplepath}
Suppose all of Assumption \ref{a:probassns} holds, and further
that the upcrossing point $\xstar$ is unique.
For any set of paths $\vrb (t)$ satisfying $\BB (\vr(T_B)) = B$
and with minimal cost,
starting from the point $\vq$, the scaled paths $\ssb (t)$ satisfy the following.
\begin{align}\label{e:TBlim}
\lim_{B\to 0} \frac{T_B}{\sqrt{B}} &= T\\
\label{e:Slim}
\lim_{B\to 0} \ssb (t) &= \vr (t),\ 0\le t\le T, \text{ uniformly in }t.
\end{align}
\end{theorem}
This theorem gives a connection between solutions to Problem A
for constant coefficient costs,
and solutions with costs depending on $\x$, but for small buffer $B$.
Proof of existence of solutions to Problem A for
non constant coefficient costs 
is the same
as the proof for Theorem \ref{bigdeal}.

{\bf Proof}:
By \cite[Prop.\ 5.46]{SWbook}, the paths
$\ssb (t)$ lie in a compact set, since they have
bounded cost and their initial points $\ssb (0)$
approach $\xstar$.
Take a convergent subsequence of $\ssb$.
By \cite[Thm.\ 5.35]{SWbook}, the costs (based on $\ell_B$)
of these paths converge to the cost of the constant coefficient path
starting at $\xstar$ and making buffer size 1 in time $T$.
But Theorem \ref{t:itsunique} shows that this path is unique
among paths that start at $\xstar$ and make buffer size 1.
Therefore any limiting sequence of the $\ssb$ must converge to the same path,
and do so uniformly, and their times must converge also.
\EndPf

We note that the point $\xstar$ will be unique if the
Markov chain can be written as a superposition
of sources, either an open model (with product-form Poisson
steady-state distribution) or a closed model (with
multinomial distribution).
In these cases, $\xstar$
is the solution of Sanov's theorem
(see \cite[Chapter 2]{SWbook}).
Explicitly, if the steady-state distribution
of the chain is $\pi(i)$, $1\le i\le K$,
then the upcrossing point is given by the
minimum entropy point (i.e., minimum Kullback-Leibler
measure) on the hyperplane $\langle\x,\va\rangle = C$.
If the level curves of this entropy
are strictly convex, then the point $\xstar$ is unique.
But it is easy to see that these level curves are strictly convex
for the two models under study, and the strict convexity
also holds in the limit as we scale $n\to\infty$.

\section{Strong uniqueness}\label{s:strong}
The preceding result may be combined with traditional variational
calculus to show strong (true) uniqueness for the solution
to Problem A from a given initial point $\vr(0)$
where $\langle\vr(0), \va\rangle = C$.
As mentioned before, this result still does not give full
uniqueness for the steady-state variational problem, because that
problem has $\vr(0) = \vq$ where $\vq$ is the unique attracting point
for the fluid limit $\z_\infty (t)$, and which satisfies
$\langle\vq,\va\rangle < C$.
Nevertheless, it is a step in that direction.

\begin{theorem}\label{t:strongunique}
If the positive cone of the 
$\{ \ei \}$
equals $\R^K$, and $\log\lai (\x)$ are Lipschitz
and $C^2$ functions, then given $\xstar$ with
$\langle\xstar,\va\rangle = C$,
there are $B_0$ and $\delta > 0$ so that for all
$B \le B_0$ and $|\x - \xstar | \le \delta$ with
$\langle\x,\va\rangle = C$,
the solution to Problem A starting at $\vr(0) = \x$ is unique.
\end{theorem}
%

%
\end{document}